\documentclass[12pt]{article} 
\usepackage{amsmath,amsfonts,times}
\advance\hoffset -0.5in

\newtheorem{theorem}{Theorem}[section]
\newtheorem{prop}[theorem]{Proposition}

\newtheorem{lemma}[theorem]{Lemma}

\newcommand{\Avec}{{\bf A}}
\newcommand{\bvec}{{\bf b}}

\newcommand{\eps}{{\varepsilon}}

\newcommand{\R}{{\mathbb R}}
\newcommand{\Uvec}{{\bf U}}
\newcommand{\uvec}{{\bf u}}
\newcommand{\vvec}{{\bf v}}
\newcommand{\wvec}{{\bf w}}
\newcommand{\xvec}{{\bf x}}
\newcommand{\Xvec}{{\bf X}}
\newcommand{\Yvec}{{\bf Y}}

\def\qed{\hbox{\hskip 6pt\vrule width6pt height7pt depth1pt
\hskip1pt}\bigskip}

\begin{document}

\title{~\\[-2cm]
ON AN ISOPERIMETRIC INEQUALITY FOR\\
A SCHR\"{O}DINGER OPERATOR DEPENDING\\
ON THE CURVATURE OF A LOOP\\[1cm]} 

\author{Almut Burchard, Lawrence E. Thomas\\[0.4cm]
University of Virginia \\ Department of Mathematics\\ Charlottesville,
Virginia 22904\\[0.5cm] 
$\{{\tt burchard,let}\}@{\tt virginia.edu}$
\date{May 6, 2005}}
\maketitle
%*****************************************************************
%*****************************************************************

\begin{abstract}
Let $\gamma$ be a smooth closed curve of length $2\pi$ in $\R^3$,
and let  $\kappa(s)$ be its curvature, regarded as a function of arc length.
We associate with this curve the one-dimensional 
Schr\"{o}dinger operator $H_\gamma=
-\frac{d^2}{ds^2} + \kappa^2(s)$ acting on the space of square
integrable $2\pi$-periodic functions.  A natural conjecture is that
the lowest spectral value $e_0(\gamma)$ of
$H_{\gamma}$ is bounded below by $1$ for any $\gamma$ 
(this value is assumed when $\gamma$ is a circle). 
We study a family of curves $\{\gamma\}$ that includes the circle 
and for which $e_0(\gamma)=1$ as well. We show that the curves 
in this family are local minimizers; i.e., $e_0(\gamma)$ can 
only increase under small perturbations leading away from the 
family.  To our knowledge, the full conjecture remains open.
\end{abstract}

%*********************************************************
\section{Introduction}
\setcounter{equation}{0}
\setcounter{theorem}{0}
\label{sec:intro}

Let $\gamma$ be a smooth closed curve of length $2\pi$
in $\R^3$, parametrized by arclength $s$.  We associate with this curve
a Schr\"odinger operator $H_\gamma$ 
on the space of square integrable, $2\pi$-periodic functions by
$$
H_\gamma \Phi(s) = -\frac{d^2\Phi(s)}{ds^2} + \kappa^2(s)\Phi(s)\ ,
$$
where $\kappa(s)$ is the curvature of $\gamma$ at $s$. Let 
\begin{equation} \label{eq:ezero}
         e_0(\gamma)= \inf\,{\rm spec}\,H_{\gamma}=\inf_{\Phi \ne 0 } 
\frac{\int_0^{2\pi} (\Phi')^2 + \kappa^2\Phi^2\, ds}
               {\int_0^{2\pi} |\Phi|^2\, ds}
\end{equation}
be the smallest eigenvalue of $H_\gamma$. It has been conjectured that
$e_0(\gamma)$ achieves its minimum
$$
    e_{\min} = \inf_{\gamma} e_0(\gamma)
$$
when $\gamma$ is a circle. In that case, $\kappa^2\equiv 1$, the 
minimizing eigenfunction $\Phi$ is constant, and 
$e_0(\gamma)=1$. But the functional assumes the same
value for an entire family ${\cal F}$ of curves
given by translations, rotations and dilations of planar
loops which
%, in the $x-y$-plane 
have tangent vector $\Uvec(s)$ proportional to 
$(\cos(s),\beta\sin(s),0)$ for some constant $\beta$ with 
$0<\beta\le 1$.  So if indeed circles are minimizers, they certainly 
are not the only minimizers.

In this article, we show that loops in the family ${\cal F}$ 
{\em locally} minimize the functional $e_0(\gamma)$ given in
Eq.~(\ref{eq:ezero}).  Small deformations about any one of these
loops cause $e_0$ to strictly increase, provided the
the loop is not simply deformed to another loop of the same family. 
This result is a first step towards understanding the landscape
in the space of curves $\{\gamma\}$ defined by the values of $e_0$.
We emphasize that the conjecture itself remains open; our results only
add credibility to it.

That $e_0(\gamma) \ge 1$ with the circle as a minimizer seems to have
been implicitly conjectured by a number of people. The conjecture was
articulated by Benguria and Loss \cite{BenLoss}, who showed
it to be equivalent to establishing the best constant for a
one-dimensional Lieb-Thirring inequality for a Schr\"{o}dinger operator
with two bound states. They did show that $e_{0}(\gamma)\geq 1/2$. We too
had made the conjecture in our work on the local existence for a
dynamical Euler elastica \cite{BT}. There, the issue of the
invertibility of $H_{\gamma}$ arises in determining the tension
of an elastic loop.  We showed that $e_0({\gamma})\geq 1/4$, 
which is in fact optimal for curves which are possibly open, and 
for which the tangent vector $\Uvec$ is $2\pi$-periodic 
and each of the components of $\Uvec$ vanishes at least once.

In related work, Harrell and Loss \cite{HL} showed that
Schr\"{o}dinger operators of the form $-\Delta - d \kappa^2$ on
$d$-dimensional hypersurfaces, with $\Delta$ the Laplace-Beltrami
operator and $\kappa$ the mean curvature, have at least two negative
eigenvalues unless the surface is a sphere (a circle in one
dimension).  Previously, Harrell ~\cite{H} had proved a similar result
for Schr\"odinger operators on embedded surfaces in ${\bf R}^3$ that
are topologically equivalent to $S^2$, with potentials given by
arbitrary definite quadratics in the principal curvatures.

Exner, Harrell, and Loss \cite{EHL} discussed a variety of isoperimetric
inequalities related to Schr\"{o}\-din\-ger operators including the
operator $H_{\gamma,g}=-d^2/ds^2+g\kappa^2(s)$ on closed curves, and
showed that, for the least eigenvalue of $H_{\gamma, g}$, the circle
is a minimizer when $g\leq 1/4$ and not a minimizer for $g>1$. 
Friedrich considered the operator with $g=1/4$ for
simple loops on the unit sphere, in connection with the Dirac operator 
on the region enclosed by such a loop~\cite{Friedrich}. 
The significance of the value $g=1$ is that  two
natural candidates for minimizing the lowest eigenvalue of $H_{\gamma,g}$
appear to exchange stability there: When $\gamma$ is 
a  circle, $\inf {\rm spec}\, H_{\gamma,g}=g$, whereas for
the extreme case of a {\em collapsed} curve ${\gamma}$, consisting 
of two straight line segments of length $\pi$ joined at their ends, 
we have $\inf {\rm spec}\, H_{\gamma,g}=1$. Such collapsed curves are
limiting points of the family ${\cal F}$.

The functional $e_0$ has no obvious convexity properties, 
and it is not amenable to standard symmetrization techniques.
One difficulty is that $\kappa^2$ cannot
be varied freely, since the condition that $\kappa$ be the curvature
of a closed curve in $\R^3$ is a complicated, nonlocal condition.
Technically, we show that the second variation of 
$e_0(\gamma_{\mu})$ is non-negative for one-parameter 
families $\gamma_{\mu}$, leading away from  a 
loop $\gamma= \gamma_{\mu}|_{\mu=0}$ in ${\cal F}$; this 
second variation is strictly positive
if the perturbation is transversal to the family.
For the case of the $\gamma$ a circle, 
where the eigenfunctions and eigenvalues of $H_{\gamma}$
are known, one can simply 
perform second order perturbation theory to show this positivity. 
For other curves in the family, the
higher eigenvalues and eigenfunctions of $H_\gamma$ 
%(apart from the ground state and its eigenvalue $1$) 
are not explicitly available, and different methods are 
needed to show the positivity.

We find it useful to rewrite the variational problem as follows.
Let $\Uvec(s)$ be the unit tangent vector to the curve, again 
parametrized by arclength $s$, let $\Phi(s)$ be the minimizing
eigenfunction, and set 
\begin{equation} \label{eq:curve-orbit}
\Xvec(s)=\Phi(s)\Uvec(s)\ ,
\end{equation}
so that $\Xvec'(s) = \Phi'(s)\Uvec(s) + \Phi(s) \Uvec'(s)$.
Since $|\Uvec(s)|\equiv1$, $\Uvec(s)\cdot \Uvec'(s)\equiv 0$,
and $|\Uvec'(s)|\equiv \kappa(s)$, we can rewrite
Eq.~(\ref{eq:ezero})
as
\begin{equation} \label{eq:def-L}
e_0(\gamma) = \frac{ \int_0^{2\pi} |\Xvec'(s)|^2\, ds}
{\int_0^{2\pi} |\Xvec(s)|^2\, ds}\ .
\end{equation}
It follows that
$$
e_{\min}= \inf \frac{ \int_0^{2\pi} |\Xvec'(s)|^2\, ds}
{\int_0^{2\pi} |\Xvec(s)|^2\, ds}\ ,
$$
where the infimum is taken over all $2 \pi$-periodic, vector-valued
functions $\Xvec$, vanishing only on a set of measure zero, with
\begin{equation} \label{eq:constraint}
\int_0^{2\pi} \frac{\Xvec(s)}{|\Xvec(s)|}\, ds = 0\ ,
\end{equation}
guaranteeing that the curve $\gamma$ with unit tangent
$\Uvec(s)=\Xvec(s)/|\Xvec(s)|$ is closed.  We will refer to 
the vector function $\Xvec(s)$ as an {\em orbit}.  Given
a vector-valued function
$\Xvec(s)$ that satisfies Eq.~(\ref{eq:constraint}),
the curve $\gamma$ can be reconstructed up to a translation
as a function $\Yvec_{\gamma}(s)\in \R^3$ by computing 
$$\Yvec_{\gamma}(s)= \int_0^s\Uvec(\tilde s)\,d\tilde s\ .$$
It is apparent that  for any choice of vectors 
$\vvec_1\ne 0$ and $\vvec_2$, the orbits 
\begin{equation} \label{eq:ellipse}
\Xvec_0(s)=  \cos(s) \vvec_1 + \sin(s) \vvec_2
\end{equation}
all satisfy the constraint in Eq.~(\ref{eq:constraint}),
and all give the same value ($e_0(\gamma)=1$) for the functional in
Eq.~(\ref{eq:def-L}). When $\vvec_1$ and $\vvec_2$ are 
linearly independent, these orbits correspond to
curves in ${\cal F}$. When $\vvec_1$ and $\vvec_2$ are linearly 
dependent, we obtain the collapsed curves
mentioned above. Our results imply the following:

\begin{theorem}\label{thm:main}
Let $\Uvec_0$ be the tangent vector to a curve $\gamma_0\in {\cal F}$,
and assume that,  for each $\mu$ sufficiently close to $0$, 
$\Uvec(\mu,s)$ describes the tangent vector of a closed curve 
of length $2\pi$  parametrized by arc length, i.e., 
$$
|\Uvec(\mu,s)|\equiv 1, \quad \int_0^{2\pi}\Uvec(\mu,s)\, ds=0\ .
$$
If $\Uvec(\mu,ds)$ has an expansion
$$
   \Uvec(\mu,s)\equiv  \Uvec_0(s)+\mu \uvec_1(s)+\mu^2\uvec_2(s)+o(\mu^2)
$$
in $H^1$, then there exists a positive number $c$ 
%depending only on $\Uvec_0$, $||\uvec_1||_{H^1}$ and $||\uvec_2||_{H^1}$ 
such that
$$
e(\gamma_{\mu})\geq e(\gamma_0)
$$
for $|\mu|<c$. The inequality is strict
unless $\gamma_\mu$ belongs again to the family ${\cal F}$.
\end{theorem}

To prove the theorem, we will show that the orbits in Eq.~(\ref{eq:ellipse})
corresponding to loops in ${\cal F}$
locally minimize the functional 
\begin{equation}\label{eq:Lagrange}
     {\cal L}(\Xvec)= 
  \frac{1}{2}\int_0^{2\pi}
\Bigl\{|\Xvec'(s)|^2-|\Xvec(s)|^2\Bigr\}\,ds
\end{equation}
subject to the constraint in Eq.~(\ref{eq:constraint}).
This implies that they locally minimize the 
functional in Eq.~(\ref{eq:def-L}).
We note in passing that the Euler-Lagrange equation for 
this minimization problem is given by
\begin{equation} \label{eq:EL}
\Xvec''(s) + \Xvec(s) = 
\frac{ |\Xvec(s)|^2\bvec - \bigl(\Xvec(s)\cdot \bvec\bigr)\Xvec(s)}
                  {|\Xvec(s)|^3} =: \Avec(s)\bvec\ ,
\end{equation}
where $\bvec\in\R^3$ is a vector of Lagrange multipliers,
and the $3\times 3$ matrix $\Avec(s)$ is
computed by differentiating the constraint
in Eq.~(\ref{eq:constraint}). 
These equations are easily seen to have first 
integrals, an {\em energy}
$$
\frac{1}{2}|\Xvec'(s)|^2 +\frac{1}{2}|\Xvec(s)|^2
- \frac{\bvec\cdot \Xvec(s)}{|\Xvec(s)|}
$$ 
and an {\em angular momentum}
$$\bvec\cdot\Xvec(s)\times\Xvec'(s)\ .
$$ We are unaware of another constant of integration which would make them an
integrable system.

In Section~\ref{sec:ellipse}, we consider 
deformations around orbits of the form
given in Eq.~(\ref{eq:ellipse}) for the generic case where $\vvec_1$ and
$\vvec_2$ are linearly independent.  These elliptical orbits
are critical points for the functional 
in Eq.~(\ref{eq:Lagrange}) even without the constraint, since they 
satisfy Eq.~(\ref{eq:EL}) with $\bvec=0$.  We show that to
second order in a parameter $\mu$ this functional can only increase for
deformations of the orbit that do not simply
transform the orbit into another elliptical orbit
new choices of $\vvec_1$ and $\vvec_2$.  The proof relies 
on an identity of elliptic integrals
which is not transparent (to us).
The section  ends with the proof of Theorem~\ref{thm:main}.

In Section~\ref{sec:collapsed}, we consider deformations about
collapsed orbits given by Eq.~(\ref{eq:ellipse})
where $\vvec_1$ is nonzero and $\vvec_2$ is a constant multiple
of $\vvec_1$. We show that the functional again increases for 
nondegenerate perturbations. Unfortunately, the analysis of 
these collapsed curves is somewhat vexing.  Their curvature
is zero along the line segments and infinite at the end points.
This forces the minimizing eigenfunctions to
vanish at these endpoints and results in a ground state of
multiplicity two so that the curve corresponds to a two-parameter
family of orbits.  We relegate the expansion of the constraint 
in Eq.~(\ref{eq:constraint}) about a collapsed critical orbit 
to the following Section~\ref{sec:expand},
the reason being that the computations are somewhat gruesome, and
their presentation would break the flow of the main arguments showing
positivity of ${\cal L}$. 

Curiously, the analysis of the second variation about the collapsed
orbits relies in part on the explicit diagonalization of the
Schr\"{o}dinger operator $K_g= -d^2/ds^2+g\sec^2(s)$, acting in
$L^2[-\pi/2,\pi/2]$ by Gegenbauer polynomials. This is discussed in the 
Appendix.

%*********************************************************
\section{Elliptical orbits}
\setcounter{equation}{0}
\setcounter{theorem}{0}
\label{sec:ellipse}

We expand  an orbit $\Xvec$  in terms of a small parameter $\mu$ as
\begin{equation} \label{eq:expand-orbit}
\Xvec(\mu,s)\equiv \Xvec_0(s)+ \mu\xvec_1(s)+\mu^2\xvec_2(s)+o(\mu^2)\, .
\end{equation}
Here, $\Xvec_0$ is a nondegenerate elliptical orbit given by
Eq.~(\ref{eq:ellipse}), $\xvec_1$ and $\xvec_2$ are
vector-valued functions in $H^1$, and the error estimate
is understood with respect to the $H^1$-norm.
Since the functional ${\cal L}$ in Eq.~(\ref{eq:def-L})
and the constraint 
in Eq.~(\ref{eq:constraint}) are symmetric under rotations, we may 
assume that
\begin{equation} \label{eq:X0}
\Xvec_0(s) = \left(\begin{array}{c}\alpha \cos(s)\\
\beta\sin{(s)}\\ 0 \end{array}\right)\ ,
\end{equation}
where $\alpha\ge \beta >0$ 
represent the major and minor semi-axes of the ellipse.
The curvature of the corresponding loop $\gamma$ is given by 
$$
\kappa(s) = \left|\frac{d}{ds} \left(
          \frac{\Xvec_0(s)}{|\Xvec_0(s)|}\right)\right|
            = \frac{\alpha\beta}{|\Xvec_0|^{2}}\ .
$$
The principal eigenvalue and eigenfunction
of the Schr\"odinger operator  $H_\gamma$ are
$$
e_0(\gamma)=1\ ,\quad \Phi(s)=|\Xvec_0(s)|  
= \sqrt{\alpha^2\cos^2(s) + \beta^2\sin^2{(s)}}\ ,
$$
and the eigenvalue-eigenvector equation reads
\begin{equation} \label{eq:EE}
H_\gamma\Phi = -\Phi''+ \frac{\alpha^2\beta^2}{\Phi^3} = \Phi\ .
\end{equation}
Expanding the functional ${\cal L}$ defined by
Eq.~(\ref{eq:Lagrange}) in powers of $\mu$,
\begin{equation}\label{eq:expand-L}
{\cal L}(\Xvec) \equiv  {\cal L}(\Xvec_0)  +\mu {\cal L}_1 + 
\mu^2{\cal L}_2 + o(\mu^2)\ ,
\end{equation}
we see that ${\cal L}_1=0$ since
$\Xvec_0$ satisfies the Euler-Lagrange equation
in Eq.~(\ref{eq:EL}).
The second variation  is given by
\begin{eqnarray*}
    {\cal L}_2
&= &\frac{1}{2}\int_0^{2\pi}\bigl\{|\xvec_1'(s)|)^2- |\xvec_1(s)|^2\bigr\}\, ds
+ \int_0^{2\pi}\bigl\{\Xvec_0'(s)\cdot \xvec_2'(s) -
                      \Xvec_0(s)\cdot \xvec_2(s)\bigr\}\, ds\\
&= & {\cal L}(\xvec_1);
\end{eqnarray*}
the contribution of $\xvec_2$ vanishes
after an integration by parts since $\Xvec_0''+\Xvec_0=0$. 
The constraint Eq.~(\ref{eq:constraint})
expanded to first order in $\mu$ implies that $\xvec_1$ satisfies the
condition
\begin{equation} \label{eq:constraint-A}
\int_0^{2\pi}\Avec(s)\xvec_1(s)\,ds = 0\ ,
\end{equation}
where
\begin{equation} \label{eq:def-Avec}
\Avec(s)= \frac{1}{|\Xvec_0|^3}\left(\begin{array}{ccc}
     \beta^2\sin^2(s) &-\alpha \beta \cos(s) \sin(s) & 0\\
     -\alpha \beta \cos(s) \sin(s) & \alpha^2 \cos^2(s)& 0\\
     0 & 0& \alpha^2\cos^2(s) + \beta^2 \sin^2(s)
\end{array}\right)
\end{equation}
is the matrix appearing in Eq.~(\ref{eq:EL}).

Consider for a moment the special case where the orbit
is a circle, $\alpha=\beta>0$. Denote the components of $\xvec_1$ by
$$
\xvec_1(s) = \left(\begin{array}{c} x_1(s)\\y_1(s)\\z_1(1)\end{array} \right)\,.
$$
The constraints in Eq.~(\ref{eq:constraint-A}) can be expressed with
the double-angle formula as
$$
\left\{
\begin{array}{l}
\displaystyle{
\int_0^{2\pi}\frac{1}{2}(1-\cos{(2s)}) x_1(s) - \frac{1}{2} \sin{(2s)}y_1(s)\, 
   ds= 0}\\
\displaystyle{\int_0^{2\pi}\frac{1}{2}(-\sin{(2s)}) x_1(s) + \frac{1}{2} 
(1+\cos{(2s)})y_1(s)\, ds= 0}\\
\displaystyle{\int_0^{2\pi}\,z_1(s) ds= 0\ .}
\end{array}\right.
$$
In other words, the zeroth and second Fourier coefficients 
of the components of $\xvec_1$ satisfy
$$
\hat x_1(\pm2) +\mp i \hat y_1(\pm2) = \hat x_1(0) + i\hat y_1(0), 
\ \hat z_1(0)=0\ .
$$
Since $\xvec_1$ is real-valued, $\hat\xvec_1(0)$ is real as well.
By the triangle inequality, $|\hat\xvec_1(0)|^2\le
2|\hat \xvec_1(\pm2)|^2$,
which  implies ${\cal L}_2\ge 0$ by
Parseval's identity.  The following proposition shows the 
corresponding statement for perturbations about general
elliptical orbits.

\begin{prop}\label{prop:localeppos}
The elliptical orbits in Eq.~(\ref{eq:X0}) locally minimize
Eq.~(\ref{eq:Lagrange}) under the constraint in Eq.~(\ref{eq:constraint})
for each $\alpha\ge \beta>0$. More precisely, 
there exists a positive constant $c=c(\alpha,\beta)$ such that
for every perturbation  $\Xvec(\mu,s)$ given by Eq.~(\ref{eq:expand-orbit})
which satisfies the constraint in Eq.~(\ref{eq:constraint}) to
order $o(\mu)$, we have 
\begin{equation} \label{claim:localeppos}
{\cal L}_2= {\cal L}(\xvec_1) \ge c(\alpha,\beta) 
\| P_{n\ne \pm 1}\xvec_1\|^2 \ , 
\end{equation}
where $P_{n\ne \pm 1}$ is the projection onto the 
space of functions whose first order Fourier coefficients vanish.
\end{prop}

\noindent {\bf Remark}: 
Variations of the form $\xvec_1(s)= 2 \mbox{\rm Re}\, e^{is} \hat \xvec(0)$
are of course along the line of critical orbits, and give
zero second variation.\\

\noindent {\sc Proof of Proposition~\ref{prop:localeppos}.}\ 
For notational convenience, we drop the subscript on $\xvec_1$
and simply write $\xvec(s)$ instead of $\xvec_1(s)$.
For the Fourier coefficients of $\xvec$ and $\Avec$, we use the convention
$$
\hat\xvec(s) = \frac{1}{\sqrt{2\pi}}\int_0^{2\pi} e^{-ins} \xvec(s)\, ds \ , 
\quad 
\hat\Avec(s) = \frac{1}{\sqrt{2\pi}}\int_0^{2\pi} e^{-ins} \Avec(s)\, ds\ .
$$
By Parseval's identity, the functional ${\cal L}$ can be expressed as 
$$
{\cal L}(\xvec) = \frac{1}{2} \sum_{n} (n^2-1)|\hat\xvec(n)|^2\ .
$$
When $\hat\xvec(0)=0$, the claim in Eq.~(\ref{claim:localeppos}) 
holds with $c=3/2$,  so we assume without
loss of generality that $\hat\xvec(0)\ne 0$.
The Fourier coefficients of $\Avec$ are nonzero only for
even $n$, since $\Avec$ is $\pi$-periodic. 
Using Parseval's identity again, we
write the constraint in Eq.~(\ref{eq:constraint-A}) as
$$
\hat\Avec(0)\hat\xvec^*(0)
   = -\sum_{n\neq 0} \hat\Avec(n)\hat\xvec^*(n)\ ,
$$
where $^*$ denotes complex conjugation.
Since the first order Fourier coefficients
of $\xvec$ contribute neither to the constraint nor
to the claim, we may assume that $\hat x(\pm 1) =0$.

The matrix $\hat\Avec(0)$ is 
invertible, since the off-diagonal elements of $\Avec(s)$
are odd in $s$ and
its diagonal elements are strictly positive, 
see Eq.~(\ref{eq:def-Avec}).
Multiplying by $\hat\Avec(0)^{-1}$ and taking the inner product
with $\hat\xvec(0)$ yields
$$
   |\hat\xvec(0)|^2= -\sum_{n\neq 0 }
   \bigl(\hat\Avec(0)^{-1}\hat\xvec(0)\bigr) \cdot \hat\Avec(n)
     \hat\xvec^*(n) = \bigl\langle -\hat\Avec^*(n)\hat\Avec(0)^{-1}\hat\xvec(0),
P\hat\xvec(n)^*\bigr\rangle_{\ell^2}\ ,
$$
where $P$ is the projection onto the nonzero Fourier modes
and $\ell^2$ denotes the space of vector-valued sequences
whose sequence of norms is square summable.
Since $\hat A(n)=0$ and $\hat\xvec(n)=0$ for $n=\pm 1$,
and $n^2-1>0$ for $n\ne 0,\pm 1$, 
we can apply the Cauchy-Schwarz inequality to obtain
\begin{equation} \label{eq:Schwarz}
      |\hat\xvec(0)|^2\leq \| (n^2-1)^{-1/2}P \hat\Avec(n)\hat\Avec(0)^{-1}
   \hat\xvec(0)\|_{\ell^2}\|(n^2-1)^{1/2}P \hat\xvec(n)\|_{\ell^2}\ .
\end{equation}
This yields the lower bound
\begin{eqnarray}\label{eq:lb1}
\nonumber
    {\cal L}(\xvec)&=& 
\frac{1}{2}\left(\| (n^2-1)^{1/2}P\hat\xvec(n)\|_{\ell^2}^2-
    |\hat\xvec(0)|^2\right)\\
\label{eq:L-eta}
&\ge&  \frac{1}{2} \left (\frac{|\hat \xvec(0)|^2}
  {\|(n^2-1)^{-1/2}P\hat\Avec(n)\hat\Avec(0)^{-1}\hat\xvec(0)\|_{\ell^2}^2} -1
\right)|\hat\xvec(0)|^2\ \\
\nonumber &\ge& \frac{\eta}{2(1-\eta)} |\hat \xvec(0)|^2\ ,
\end{eqnarray}
where  $\eta$ 
is the lowest eigenvalue of the $3\times 3$ matrix
\begin{equation} \label{eq:dual}
D= \hat \Avec(0)^{-1} \Bigl\{\sum_{n \ne  \pm 1} \frac{1}{1-n^2}
\hat\Avec (n) \hat \Avec(n)^*\Bigr\}\hat\Avec(0)^{-1}\ .
\end{equation}
Note that the idenitity matrix is included
as the $n=0$ term in the definition of
$D$. Clearly $\eta<1$ since $D$ is the 
identity minus a positive
definite matrix. We will show that
$\eta> 0$ by verifying that the sum inside the braces of 
Eq.~(\ref{eq:dual}) is a positive definite matrix.

We express this sum as a convolution integral.  In order 
to invert the Fourier multiplication
operator $1-n^2$ on the space of functions whose odd Fourier
coefficients vanish, we need to solve the equation 
$$y''+y=f$$
 on the space of $\pi$-periodic functions.  Since $K(s) =
\frac{1}{4}|\sin{(s)}|$ satisfies $K''(s) + K(s) = \frac{1}{2}
(\delta_0+\delta_\pi)$, the unique $\pi$-periodic solution is given by
$$ 
K*f(s) = \int_0^{2\pi} K(s-t)f(t)\, dt\ ,
$$
and so
\begin{equation}\label{eq:brace}
\sum_{n \ne  \pm 1} \frac{1}{1-n^2}
\hat\Avec (n) \hat \Avec(n)^*= \frac{1}{4}
\int_0^{2\pi}\int_0^{2\pi} \Avec(s)\Avec(t)|\sin(s-t)|\, 
ds dt\ .
\end{equation}
From the expression for $\Avec(s)$ in Eq.~(\ref{eq:def-Avec}) it is apparent that the off-diagonal terms in
$\Avec(s)\Avec(t)$ change sign if $(s,t)$ is replaced by $(-s,-t)$
and hence integrate to zero. Thus the 
expression in Eq.~(\ref{eq:brace}) is actually diagonal with diagonal 
entries given by 
\begin{eqnarray} 
\nonumber
I_1 &=& \Bigl\langle A_{11},K*A_{11}\Bigr\rangle_{L^2} 
      + \Bigl\langle A_{12},K*A_{12}\Bigr\rangle_{L^2}\\
\label{eq:I123}
I_2 &=& \Bigl\langle A_{22},K*A_{22}\Bigr\rangle_{L^2} 
      + \Bigl\langle A_{12},K*A_{12}\Bigr\rangle_{L^2}\\
\nonumber
I_3 &=& \Bigl\langle A_{33},K*A_{33}\Bigr\rangle_{L^2}\ ,
\end{eqnarray}
where $A_{ij}$ is the $ij$-th entry of $\Avec$.
It just remains to show positivity of these $I_j$'s. Clearly,
$$
I_3=  \frac{1}{4}\int_0^{2\pi}\int_{0}^{2\pi} 
|\Xvec_0(s)|^{-1}|\Xvec_0(t)|^{-1}|\sin{(s-t)}|\, ds dt>0\ ,
$$
and we note that
$$
  I_1+I_2= \int_0^{2\pi}\int_0^{2\pi}\frac{(\alpha^2\cos{(s)}\cos{(t)}+
\beta^2\sin{(s)}\sin{(t)})^2}
{|\Xvec_0(s)|^3|\Xvec_0(t)|^3} |\sin(t-s)|\,ds dt>0
$$
since the integrands are nonnegative.
It follows from Lemma~\ref{lem:amazing}, which is proved below,
that $I_1= \frac{\beta^2}{\alpha^2+\beta^2}(I_1+I_2)$ and
$I_2= \frac{\alpha^2}{\alpha^2+\beta^2}(I_1+I_2)$ are both
positive.  Since $\hat\Avec(0)$ is a diagonal matrix with positive entries,
we conclude from Eq.~(\ref{eq:dual}) that $\eta>0$, and 
hence ${\cal L}_2>0$.  \hfill\qed

In the proof of Proposition~\ref{prop:localeppos}, we used that
$I_1$ and $I_2$ are positive multiples of $I_1+I_2$.
This is a consequence of the following identity which we state 
as a lemma. We have no geometric insight why this identity 
should hold; it was discovered numerically.

\begin{lemma} \label{lem:amazing} 
The integrals in Eq.~(\ref{eq:I123}) satisfy $\alpha^2 I_1=\beta^2I_2$.
\end{lemma}

\noindent {\sc Proof.}\ 
The lemma clearly holds for $\alpha=\beta>0$, since then
$I_2$ can be obtained from $I_1$ by replacing $(s,t)$ 
with $(s+\pi/2,t+\pi/2)$.  For $\alpha>\beta>0$, we write
$$
A_{11}(s)  
=\frac{\beta^2\sin^2{(s)}}{|\Xvec_0(s)|^3}
= -\frac{\beta^2}{(\alpha^2-\beta^2)} |\Xvec_0(s)| ^{-1}
+ \frac{\alpha^2\beta^2}{(\alpha^2-\beta^2)}|\Xvec_0(s)|^{-3}\ .
$$
Since $\alpha^2\beta^2 K*|\Xvec_0(s)|^{-3}=|\Xvec_0(s)|$
by Eq.~(\ref{eq:EE}) and the definition of $K$, we have
\begin{equation} \label{eq:I1-a}
\begin{array}{lcl}
\displaystyle{ \Bigl\langle A_{11}, K*A_{11}\Bigr\rangle_{L^2}} &= &
\displaystyle{
\frac{\beta^4}{(\alpha^2-\beta^2)^2 }
     \Bigl\langle |\Xvec_0|^{-1},K*|\Xvec_0|^{-1}\Bigr\rangle_{L^2}}\\
&& \hskip -2cm - \ 
\displaystyle{ 2 \frac{\beta^2}{(\alpha^2-\beta^2)^2} 
     \Bigl\langle |\Xvec_0|^{-1}, |\Xvec_0|\Bigr\rangle_{L^2}
+ \frac{\alpha^2\beta^2}{(\alpha^2-\beta^2)^2} 
     \Bigl\langle |\Xvec_0|^{-3},|\Xvec_0|\Bigr\rangle_{L^2} \ .}
\end{array}
\end{equation}
For the second term in $I_1$, we compute
$$
A_{12}(s) = -\frac{\alpha\beta \cos{(s)}\sin{(s)}}{|\Xvec_0(s)|^3} = 
-\frac{\alpha\beta}{\alpha^2-\beta^2} \frac{d}{ds} |\Xvec_0(s)|^{-1}\ ,
$$ 
which gives 
$$
\frac{d}{ds} K*A_{12}
= -\frac{\alpha\beta }{\alpha^2-\beta^2} \frac{d^2}{ds^2} K*|\Xvec_0|^{-1}
= \frac{\alpha\beta}{\alpha^2-\beta^2} K*|\Xvec_0|^{-1}
- \frac{\alpha\beta }{\alpha^2-\beta^2} |\Xvec_0|^{-1} \ .
$$
by the definition of $K$. With an integration by parts, we see that
\begin{eqnarray}\label{eq:I1-b}
\lefteqn{\Bigl\langle A_{12}, K*A_{12}\Bigr\rangle_{L^2}}\\
  &=&  \frac{\alpha^2\beta^2}{(\alpha^2-\beta^2)^2}
         \Bigl\langle|\Xvec_0|^{-1}, K*|\Xvec_0|^{-1}\Bigr\rangle_{L^2} 
- \frac{\alpha^2\beta^2}{(\alpha^2-\beta^2)^2} 
         \Bigl\langle|\Xvec_0|^{-1}, |\Xvec_0|^{-1}\Bigr\rangle_{L^2}\ .
\nonumber
\end{eqnarray}
Adding Eqs.~(\ref{eq:I1-a}) and~(\ref{eq:I1-b}), we obtain
$$
I_1= \beta^2 \biggl\{ \frac{\alpha^2+\beta^2}{(\alpha^2-\beta^2)^2} 
         \Bigl\langle|\Xvec_0|^{-1} , K*|\Xvec_0|^{-1}\Bigr\rangle_{L^2}
- 2 \frac{1}{(\alpha^2-\beta^2)^2}\biggr\}
$$
In the same way, we compute
$$
I_2= \alpha^2 \biggl\{\frac{\alpha^2+\beta^2}{(\alpha^2-\beta^2)^2} 
         \Bigl\langle|\Xvec_0|^{-1}, K*|\Xvec_0|^{-1}\Bigr\rangle_{L^2}
- 2 \frac{1}{(\alpha^2-\beta^2)^2}\biggr\} \ ,
$$
which proves the lemma.\hfill\qed

\bigskip
The lower bound on ${\cal L}_2$ in Proposition~\ref{prop:localeppos} 
deteriorates when the elliptical orbit $\Xvec_0$ collapses.
Fix $\alpha=1$, and let $\beta\to 0$.
By an analysis of the integrands in Eq.~(\ref{eq:I123}),
particularly near $s,t= \pm \pi/2$, we find that 
\begin{eqnarray*}
       I_1&\sim& \beta^2\ln(1/\beta),\\ 
       I_2&\sim& \ln(1/\beta),\\ 
       I_3&\sim& \ln(1/\beta), 
\end{eqnarray*}
and similarly
$$
\hat\Avec(0) \sim
\left(\begin{array}{ccc}
       1 & 0 & \\ 
       0 & \ln(1/\beta) & 0 \\
       0 & 0 & \ln(1/\beta)
\end{array}\right)
$$ 
It follows that the lowest eigenvalue of the diagonal matrix $D$ in
Eq.~(\ref{eq:dual}) is given by the entry involving $I_1$, and so, by
Eq.~(\ref{eq:L-eta}),
$$
     {\cal L}_2 \ge \frac{\eta}{2(1-\eta)}
\sim \beta^2\ln(1/\beta) |\hat \xvec_1(0)|^2\ .
$$
On the other hand,
$${\cal L}_2 \ge \frac{3}{2}
 ||P_{n\ne 0,\pm 1}\xvec_1||^2- \frac{1}{2} |\hat \xvec_1(0)|^2\ ,
$$
using the first line of Eq.~(\ref{eq:lb1}). 
Interpolating between these two inequalities
we obtain
$$
{\cal L}_2 \ge c \beta^2 \ln(1/\beta)\, { || P_{n\ne \pm 1}\xvec_1||^2} 
\qquad 
(\mbox{$\alpha=1$, $\beta\to 0$})\ ,
$$
where $c$ is an absolute constant.  
Since Eq.~(\ref{eq:Schwarz}) can hold with equality, 
the lowest eigenvalue of ${\cal L}$ on the space of 
functions whose first order Fourier coefficients vanish is 
also bounded above by a constant multiple of 
$\beta^2 \ln(1/\beta)$.

\bigskip\noindent{\sc Proof of Theorem~\ref{thm:main}.} \ 
Let $\Uvec(\mu,s)$ be as in the statement of
the theorem, and let  $\Phi(\mu,s)$ be the normalized
minimizing eigenfunction for the corresponding curve
$\gamma(\mu)$. Since the ground state of $H_\gamma$
is simple, we may expand $\Phi(\mu,s)$ in $H^1$ as
$$
\Phi(\mu,s)\equiv \Phi_0(s) + \mu \phi_1(s) + \mu^2\phi_2(s) + o(\mu^2)\ .
$$
The corresponding orbit is given by
$\Xvec(\mu,s)=\Phi(\mu,s)\Uvec(\mu,s)$, see
Eq.~(\ref{eq:curve-orbit}), which has an expansion
as in Eq.~(\ref{eq:expand-orbit}) with
\begin{eqnarray*}
\Xvec_0(s) &=& \Phi_0(s) \Uvec_0(s)\\
\xvec_1(s) &=& \phi_1(s) \Uvec_0(s) + \Phi_0(s)\uvec_1(s)\\
\xvec_2(s) &=& \phi_2(s) \Uvec_0(s) + \phi_1(s)\uvec_1(s)
             + \Phi_0(s)\uvec_2(s)\ .
\end{eqnarray*}
Since the unperturbed curve
$\Uvec_0$ belongs to the family ${\cal F}$, 
we may assume by performing a suitable rotation and
translation that $\Xvec_0(s)$
satisfies Eq.~(\ref{eq:X0}).  By
Proposition~\ref{prop:localeppos}), there exists
a constant $c>0$ such that
${\cal L}(\Xvec(\mu,s)) \ge 0 $ for $|\mu|<c$, with strict inequality
if the variation is transversal to the family 
${\cal F}$. The claim now follows from the definition of $L$ in
Eq.~(\ref{eq:def-L}).
\hfill\qed

%************************************************
\section{Collapsed orbits}
\setcounter{equation}{0}
\setcounter{theorem}{0}
\label{sec:collapsed}

If the vectors $\vvec_1$ and $\vvec_2$ defining the elliptical
orbits in Eq.~(\ref{eq:ellipse}) are linearly dependent, then the
corresponding curve collapses into a pair of straight line segments
joined at the ends. The associated Schr\"odinger operator is just the
second derivative operator acting on $2\pi$-periodic functions in $H^1$
which vanish at $\pi/2$ and $3\pi/2$.  The lowest eigenvalue of this
operator is $e_0=1$, and has multiplicity two, and the eigenfunctions
are multiples of
$$
\cos_{\alpha\beta}(s) = \left\{ \begin{array}{ll}
            \alpha \cos{(s)},\quad &\mbox{if}\  -\pi/2 \le s\le \pi/2\\
            \beta \cos{(s)}, & \pi/2\le s\le 3\pi/2\ ,
\end{array}\ \right.
$$
where $\alpha$ and $\beta$  are constants.  The corresponding orbits 
are given by
\begin{equation} \label{eq:cos-rho}
\Xvec_0(s) = \left(\begin{array}{c} \cos_{\alpha\beta}{(s)} \\ 0 \\ 0 
\end{array}\right)\ .
 \end{equation}
In this section, we show that these collapsed orbits also locally
minimize the functional ${\cal L}$.  We consider perturbations
around an orbit $\Xvec_0$  given by Eq.~(\ref{eq:cos-rho})
with $\alpha>0$ and $0^+\le\beta\le\alpha$.  We expand the 
perturbation to order $o(\mu^2)$ in $H^1$ as in 
Eq.~(\ref{eq:expand-orbit}).  Expanding ${\cal L}$ as 
in Eq.~(\ref{eq:expand-L}), we obtain for the first variation
\begin{eqnarray}\label{eq:L1.1}
{\cal L}_1 &=& \int_0^{2\pi} 
\left\{ \Xvec_0'(s)\cdot \xvec_1'(s) -
\Xvec_0(s)\cdot \xvec_1(s)\right\}\, ds \nonumber\\
&=& -(\alpha-\beta) \bigl\{ x_1(\pi/2) +x_1(3\pi/2) \bigr\}\ .
\end{eqnarray}
We have integrated by parts on each of the intervals $[-\pi/2,\pi/2]$
and $[\pi/2,3\pi/2]$ and used that $\Xvec_0''+\Xvec_0=0$ in the
interior of these intervals.  Note that ${\cal L}_1$
vanishes when $\alpha=\beta$.  For $\alpha\ne \beta$ the boundary terms 
can be of either sign, indicating that these orbits are 
not critical for ${\cal L}$ without constraints. 
We will show that ${\cal L}$ can only increase under
small non-degenerate deformations that respect the constraint
in Eq.~(\ref{eq:constraint}).

\begin{prop} \label{prop:L1} Let 
$\Xvec_0$ be an orbit defined by by Eq.~(\ref{eq:cos-rho}) with
$\alpha>0$ and $0^+\le\beta\le\alpha$.  Consider perturbations of
$\Xvec_0$ given by
$$
\Xvec(\mu,s) \equiv  \Xvec_0(s) + \mu\Xvec_1(s) + o(\mu)
$$
in $H^1$, and let the corresponding expansion of ${\cal L}$ be given by
$$
{\cal L}(\Xvec) = {\cal L}(\Xvec_0) + \mu {\cal L}_1 + o(\mu)\ .
$$
If the first component of the constraint in Eq.~(\ref{eq:constraint}) 
is satisfied to order $o(\mu)$, then
$\mu {\cal L}_1\ge 0$. It is strictly positive
unless either $\alpha=\beta>0$ or 
$\xvec_1(\pi/2)=\xvec_1(3\pi/2)=0$.  
\end{prop}

\noindent{\sc Proof.} \ 
As mentioned in the introduction, we will
need an expansion of the constraint 
in Eq.~(\ref{eq:constraint}). This expansion
is provided by Lemma~\ref{lem:1storder} in the next section. 

Consider the first case where $\alpha>0$ and $\beta=0^+$.  
Denote the components of the perturbed orbit by
\begin{equation} \label{eq:components-1}
\Xvec(\mu,s) = \left(\begin{array}{c} X(\mu,s)\\Y(\mu,s)\\Z(\mu,s)
\end{array}\right)\ ,\quad 
\xvec_1(s)=\left(\begin{array}{c} x_1(s)\\y_1(s)\\z_1(s)
\end{array}\right)\ .
\end{equation}
By Lemma~\ref{lem:1storder},
the contribution of the interval $[-\pi/2,\pi/2]$ to the
first component of the integral in Eq.~(\ref{eq:constraint}) has an expansion
\begin{equation} \label{eq:1storderx-a}
\int_{-\pi/2}^{\pi/2}\frac{X(\mu,s)}
{|\Xvec(\mu,s)|}\, ds
= \pi+  O(\mu)\ .
\end{equation}
The contribution of $[\pi/2, 3\pi/2]$ is given by
\begin{equation} \label{eq:1storderx-b}
\int_{\pi/2}^{3\pi/2}
\frac{X(\mu,s)}{|\Xvec(\mu,s)|}\, ds = 
\int_{\pi/2}^{3\pi/2}
\frac{\mu x_1(s)+o(\mu)}{|\mu\xvec_1(s)+o(\mu)|}\,ds \ge -\pi\ .
\end{equation}
If $\mu x_1(\pi/2)>0$, then $X(\mu,s)$ is greater than zero on a set
whose measure does not go to zero as $\mu\to 0$. The same is true if
$x_1(3\pi/2)>0$.  Similarly, if $y_1$ or $z_1$ is nonzero for some
$s\in [\pi/2,3\pi/2]$, then by the continuity of these functions, the
integrand differs from $-1$ by at least some fixed positive value on a
set whose measure does not go to zero as $\mu\to 0$.  In either case,
the integral then would strictly exceed $-\pi+\eps$ for some $\eps>0$
for all sufficiently small values of $\mu$.  Adding
Eqs.~(\ref{eq:1storderx-a}) and~(\ref{eq:1storderx-b}), we see that if
the constraint in Eq.~(\ref{eq:constraint}) is satisfied to order
$\mu$, then $\mu x_1(\pi/2)\le 0$, $\mu x_1(3\pi/2)\le 0$, and $y_1$
and $z_1$ vanish identically on $[\pi/2,3\pi/2]$.  The claim follows
now directly from the expression for ${\cal L}_1$ in
Eq.~(\ref{eq:L1.1}).

If $\beta>0$, we use Lemma~\ref{lem:1storder}
to expand the integral in  Eq.~(\ref{eq:constraint}) over $[\pi/2,3\pi/2]$
as well as $[-\pi/2,\pi/2]$, 
\begin{eqnarray} \label{eq:1storderx-d}
\int_{0}^{2\pi}\frac{X(\mu,s)}
{|\Xvec(\mu,s)|}\, ds 
&=&    \mu\Bigl(\frac{1}{\alpha}+\frac{1}{\beta}\Bigr) 
\Bigl(x_1(-\pi/2)+ x_1(\pi/2)\Bigr) \\
&& \qquad
- |\mu|\Bigl(\frac{1}{\alpha}-\frac{1}{\beta}\Bigr) 
\Bigl(|\xvec_1(\pi/2)|+ |\xvec_1(3\pi/2)|\Bigr) 
+ o(\mu) \,.
\nonumber
\end{eqnarray}
Setting the leading term in Eq.~(\ref{eq:1storderx-d}) equal to zero,
solving for $x_1(\pi/2)+x_1(3\pi/2)$ and inserting the result into
Eq.~(\ref{eq:L1.1}), we see that
$$ 
\mu{\cal L}_1 
=|\mu|\frac{ (\alpha-\beta)^2}{\alpha+\beta} \bigl( |\xvec_1(\pi/2)| +      
             |\xvec_1(3\pi/2)|\bigr) 
\ge 0\ ,
$$
as claimed.  \hfill\qed

\bigskip\noindent 
If $\alpha=\beta$ or $\xvec_1(\pi/2)=\xvec_1(3\pi/2)=0$, we must 
work to higher order in $\mu$ to detect positivity of ${\cal L}$. 
Expanding ${\cal L}$
to second order in $\mu$ yields with a similar computation
as in Eq.~(\ref{eq:L1.1})
\begin{equation} \label{eq:L2.1}
{\cal L}_2 = \frac{1}{2}\int_0^{2\pi}
\left\{ |\xvec'_1(s)|^2 - |\xvec_1(s)|^2\right\}\, ds 
- (\alpha-\beta) \bigl\{x_2(\pi/2) +x_2(3\pi/2) \bigr\} \ .
\end{equation}
Our next result is that the second variation of the functional
is nonnegative whenever the first variation vanishes.

\begin{prop} \label{prop:L2}
Let $\Xvec_0$ be given by Eq.~(\ref{eq:cos-rho}),
and let $\Xvec(\mu,s)$ be an $H^1$-perturbation
of $\Xvec_0$, given by an expansion
as in Eq.~(\ref{eq:expand-orbit}).
Assume that the first component of the constraint in Eq.~(\ref{eq:constraint}) 
is satisfied to order $o(\mu^2)$, and the second and third
components of Eq.~(\ref{eq:constraint})
are satisfied to order $o(\mu)$.
Consider the corresponding expansion of ${\cal L}$ given by
Eq.~(\ref{eq:expand-L}).
If ${\cal L}_1=0$, then ${\cal L}_2\ge 0$. If the perturbation
is transversal to the family of collapsed orbits, then
${\cal L}_2>0$.

\end{prop}

\noindent{\sc Proof.} \ 
Let $\alpha\ge\beta\ge 0^+$, $\Xvec$, and $\Xvec_0$ be as in the
statement of the theorem. Denote the components of
the vector-valued functions appearing in the
Eq.~(\ref{eq:expand-orbit}) by
\begin{equation}\label{eq:components-2}
\Xvec(\mu,s) = \left(\begin{array}{c} X(\mu,s)\\Y(\mu,s)\\Z(\mu,s)
\end{array}\right)\ ,\quad
\xvec_1(s) = \left(\begin{array}{c} x_1(s)\\y_1(s)\\z_1(s)
\end{array}\right)\ ,\quad
\xvec_2(s) = \left(\begin{array}{c} x_2(s)\\y_2(s)\\z_2(s)
\end{array}\right)\  .
\end{equation}
Since ${\cal L}_1=0$, we have by
Proposition~\ref{prop:L1} that either $\alpha=\beta$ or
$\xvec_1(\pi/2)=\xvec_1(3\pi/2)=0$.
When $\alpha=\beta$, we invoke the first component of
the constraint to order $o(\mu)$ and the second and third components
to order $o(1)$ and use Lemma~\ref{lem:1storder})
to conclude that $\xvec_1(\pi/2)=\xvec_1(3\pi/2)=0$ as well.
In either case, the integral involving $x_1$ in 
Eq.~(\ref{eq:L2.1}) is strictly
positive, unless the restrictions of $\xvec_1$ to $[-\pi/2,\pi/2]$ 
and $[\pi/2,3\pi/2]$ are multiples of $\cos (s)$. 
Expanding the second and third component of the
constraint in Eq.~(\ref{eq:constraint}) to order $o(1)$
and using Lemma~\ref{lem:1storder}, we see that
then $y_1$ and $z_1$ are multiples of $\cos_{\alpha\beta}$,
i.e., the variation is in the direction of the family of collapsed orbits.
When $\alpha=\beta>0$, this concludes the argument.
For $\alpha>\beta$, the terms containing $y_1$ and $z_1$ will 
be used to balance the terms containing~$\xvec_2$.

Consider first the case where $\alpha>0$ and $\beta=0$.
By Lemma~\ref{lem:2ndorder},  the contribution 
of the interval $[-\pi/2,\pi/2]$ to the integral in 
Eq.~(\ref{eq:constraint}) satisfies
\begin{equation} \label{eq:2ndx-a} 
\int_{-\pi/2}^{\pi/2}\frac{X(\mu,s)}
{|\Xvec(\mu,s)|}\, ds
=    \pi + O(\mu^2)
\end{equation}
If $x_2(\pi/2)>0$, then it follows
from the continuity estimate in Eq.~(\ref{eq:Holder})
that $X(\mu,s)=\mu x_1(s)+\mu^2
x_2(s)+o(\mu^2)$ is nonnegative on an interval 
$[\pi/2,s^*(\mu)]$,
where $s^*(\mu)-\pi/2=\mu^2/o(1)$ as $\mu\to 0$.
It follows that the contribution
of the interval $[\pi/2,3\pi/2]$ satisfies 
\begin{equation} \label{eq:2ndx-b}
\int_{\pi/2}^{3\pi/2} \frac{X(\mu,s)}{|\Xvec(\mu,s)|}\, ds
= \int_{\pi/2}^{3\pi/2} \frac{ \mu x_1(s)+\mu^2 x_2(s) + o(\mu^2)}
{ |\mu\xvec_1(s)+\mu^2\xvec_2(s) + o(\mu^2)|}
\ge -\pi + \frac{\mu^2}{o(1)}\ .
\end{equation}
Adding Eqs.~(\ref{eq:2ndx-a}) and~(\ref{eq:2ndx-b}), we see that
then the constraint in Eq.~(\ref{eq:constraint})
cannot be satisfied to order $o(\mu^2)$.
Therefore $x_2(\pi/2)$ and similarly $x_2(3\pi/2)$ cannot be positive. 
The claim now follows directly from Eq.~(\ref{eq:L2.1}).

When $\alpha \ge  \beta>0$, we use Lemma~\ref{lem:2ndorder} to expand
the first component of the constraint in
Eq.~(\ref{eq:constraint})
over the entire interval $[0,2\pi]$,
\begin{eqnarray}\label{eq:2ndx-c} 
\int_{0}^{2\pi}\frac{X(\mu,s)} {|\Xvec(\mu,s)|}\, ds
&=&  \mu^2 \Bigl\{-\frac{1}{2}\int_{0}^{2\pi}
\frac{\mbox{\rm sign}\,(\cos{(s)}) }{|\cos^2_{\alpha\beta}(s)|}
\bigl(y_1^2(s)+z_1^2(s)\bigr)\,ds \\
&& \hskip -3cm+ 
\Bigl(\frac{1}{\alpha}+\frac{1}{\beta}\Bigr) 
\Bigl(x_2(\pi/2)+ x_2(3\pi/2)\Bigr) 
- \Bigl(\frac{1}{\alpha}-\frac{1}{\beta}\Bigr) 
\Bigl(|\xvec_2(\pi/2)|+ |\xvec_2(3\pi/2)|\Bigr)
\Bigr\} +o(\mu^2)\,\ .\nonumber
\end{eqnarray}
The integral on the right hand side is well-defined
by Lemma~\ref{lem:lowerbound} of the Appendix.  
To enforce the constraint in Eq.~(\ref{eq:constraint}),
we set the leading term in Eq.~(\ref{eq:2ndx-c}) equal to zero and
solve for $x_2(\pi/2)+x_2(3\pi/2)$.
Inserting the resulting expression into Eq.~(\ref{eq:L2.1}) yields
\begin{eqnarray}
\nonumber
{\cal L}_2 
&=& 
\frac{1}{2} \int_0^{2\pi} 
\left\{ |\xvec_1'(s)|^2-|\xvec_1(s)|^2
+ g_{\alpha\beta}(s) \sec^2(s)\bigl(y_1^2(s)+z_1^2(s)\bigr) \right\}\, ds\ ,\\
&&
\label{eq:L2b} 
+\frac{(\alpha-\beta)^2}{\alpha+\beta}
\bigl\{|\xvec_2(\pi/2)| + |\xvec_2(3\pi/2)|\bigr\}\\
&&\nonumber \quad
\end{eqnarray} 
where 
\begin{equation} \label{eq:g-alphabeta}
g_{\alpha\beta}(s)\equiv \left\{\begin{array}{ll}
          - \frac{\beta(\alpha-\beta)}{\alpha (\alpha+\beta)}\quad 
    & -\pi/2\leq s < \pi/2\\
          \ \ \frac{\alpha(\alpha-\beta)}{\beta(\alpha+\beta)}& 
\phantom{-}\pi/2\leq s<3\pi/2\ .
        \end{array}\right.
\end{equation}
The terms involving $\xvec_2$ in Eq.~(\ref{eq:L2b}) 
are clearly nonnegative.
The part of the integral involving the first component
$x_1$ is nonnegative because $x_1$ 
vanishes at $\pi/2$ and $3\pi/2$. 

\medskip 
To analyze the contribution
of $y_1$ to the integral in Eq.~(\ref{eq:L2b}), we invoke
the second component
of the constraint in Eq.~(\ref{eq:constraint}) to order
$o(\mu)$. By Lemma~\ref{lem:2ndorder}, 
$$
\int_0^{2\pi}\frac{Y(\mu,s)}{|\Xvec(\mu,s)|} \,ds
= \mu\int_{0}^{2\pi}\frac{y_1(s)}{|\cos_{\alpha\beta}(s)|}\,ds+o(\mu)\ .
$$
The corresponding statements hold for the third component, $z_1$.
Thus, we minimize
\begin{equation} \label{eq:secsqrpot}
  \int_{0}^{2\pi}\left\{\left(w'(s)\right)^2+
g_{\alpha\beta}(s) \sec^2(s) w^2(s)\right\}\,ds
\end{equation}
on the space of $2\pi$-periodic functions in $H^1$-functions that vanish
at $\pi/2$ and $3\pi/2$  
subject to the constraints that
 \begin{equation}\label{eq:secsqrpot-constraint}
 ||w||_2^2=1\ ,\quad
    \int_0^{2\pi}\frac{w(s)}{|\cos_{\alpha\beta}(s)|}\,ds = 0\ .
\end{equation}
We will prove that the minimum is $1$, thereby showing that
the total contributions of $y_1$ and $z_1$ to Eq.~(\ref{eq:L2b}) 
are nonnegative.

The Euler-Lagrange equation for the minimization problem in
Eqs.~(\ref{eq:secsqrpot})-(\ref{eq:secsqrpot-constraint}) is
given by
\begin{equation} \label{eq:eta-nu}
Kw(s) :=  -\frac{d^2w(s)}{ds^2} + g_{\alpha\beta}(s)\sec^2 (s) w (s)
 = \frac{\nu}{|\cos_{\alpha\beta}(s)|} +  \eta w\ ,
\end{equation}
where $\eta = \bigl(w,Kw\bigr)$  is the value 
of the functional, and $\nu$ is a Lagrange multiplier. 
We verify by direct computation that 
$$
    w_0(s)= -\frac{\nu\alpha(\alpha+\beta)}{\beta(\alpha-\beta)}
\cos_{\alpha\beta}(s)
$$
solves Eq.~(\ref{eq:eta-nu}) with $\eta=1$. This shows that
$\eta=1$ is a critical value of the functional.

Since $g_{\alpha\beta}>-1/4$ by Eq.~(\ref{eq:g-alphabeta}), we can apply
Lemma~\ref{lem:lowerbound} from the appendix to see that
the operator $K$ is bounded below and has compact resolvent.
The spectrum of $K$ consists of an increasing
sequence of eigenvalues $\lambda_0,\lambda_1,\dots$
with $\lambda_n\to\infty$.  The spectrum of $K$ is the union of the 
spectra of its restrictions to $[-\pi/2,\pi/2]$ and $[\pi/2,3\pi/2]$, 
which are determined explicitly in the appendix. 
It follows from Eq.~(\ref{eq:lambda01}) that 
$\lambda_0>1/4$ and $\lambda_1>1$.

Furthermore, a solution of the minimization problem in 
Eqs.~(\ref{eq:secsqrpot})-(\ref{eq:secsqrpot-constraint})
exists. In fact, the  constrained functional
has an infinite sequence of critical values $\eta_0\le\eta_1\le \dots$,
for which the Euler-Lagrange equation in Eq.~(\ref{eq:eta-nu})
has a nontrivial solution. If $P$ is the projection onto 
the orthogonal complement of $1/|\cos_{\alpha\beta}|$ in $L^2$, then
these critical values are just the eigenvalues of
the operator $PKP$. By the minimax characterization of 
eigenvalues of self-adjoint operators, 
the second-lowest critical value $\eta_1$ satisfies
\begin{eqnarray*}
\eta_1
&\geq &\min_{\{D: D\perp 1/|\cos_{\alpha\beta}|\}} \; 
       \max_{\{w\in D:\|w\|=1\}}\;
                              \bigl\langle w, Kw\bigr\rangle_{L^2} \\
&\geq &\min_{\{D\}}\quad \max_{\{w\in D: \|w\|=1\}}\;
                              \bigl\langle w, Kw\bigr\rangle_{L^2} \\
&=& \lambda_1 >1 \ .
\end{eqnarray*}
Here, $D$ runs over two-dimensional subspaces of $L^2$,
see Theorem 12.1 of~\cite{LL}, Eq. (5). 

We conclude that $w_0$ is indeed the minimizer, 
and $\eta_0= 1$ is the minimum value.  Since
$\eta_1$ can also be characterized by $$\eta_1 = 
\min\bigl\{ (w,Kw):\ ||w||^2=1, w\perp w_0, w\perp 1/|\cos_{\alpha\beta}
\bigr\}\ ,
$$
the functional in Eq.~(\ref{eq:secsqrpot}) is bounded
below on the subspace of functions perpendicular 
to $1/|\cos_{\alpha\beta|}$ by 
$$
\bigl\langle w,Kw\bigr\rangle _{L^2} \ge ||w||_{L^2}^2 + (\eta_1-1) 
    \bigl\{ ||P_{w_0^\perp}w||^2\bigr\}\ ,
$$
where $P_{w_0^\perp}$ is the projection onto the subspace
orthogonal to $w_0$.
\hfill\qed

%***********************************
\section{The constraint integrals near a collapsed orbit}
\setcounter{equation}{0}
\setcounter{theorem}{0}
\label{sec:expand}

In this section we consider two expansions for $\Xvec(\mu,s)$
about a singular orbit $\Xvec_0$, as given
in Eq.~(\ref{eq:cos-rho}).  The calculations are 
summarized in the following two lemmas.

\begin{lemma}\label{lem:1storder}
Assume that  a vector-valued function $\Xvec$ on the interval
$[-\pi/2,\pi/2]$ satisfies
\begin{equation}\label{eq:xvecform1}
      \Xvec(\mu,s)=
\left(\begin{array}{c}\alpha\cos{(s)}\\ 0\\0\end{array}\right)
+\mu\xvec_1(s)+o(\mu)\  ,
\end{equation}
in $H^1$.  
Then, using the notation of Eq.~(\ref{eq:components-1}),
\begin{eqnarray} \label{eq:1storderx}
\int_{-\pi/2}^{\pi/2}\frac{X(\mu,s)}
{|\Xvec(\mu,s)|}\, ds
&=& \pi\,  \mbox{\rm sign}({\alpha}) + 
\frac{\mu}{|\alpha|} \bigl(x_1(-\pi/2)+ x_1(\pi/2)\bigr) \\
&&\quad - \frac{|\mu|}{\alpha}
\bigl(|\xvec_1(-\pi/2)|+ |\xvec_1(\pi/2)|\Bigr) 
+ o(\mu)
\nonumber
\end{eqnarray}
and 
\begin{eqnarray} \label{eq:1stordery}
\int_{-\pi/2}^{\pi/2}\frac{1}{|\Xvec(\mu,s)|}\left(\!\begin{array}{c}
Y(\mu,s)\\ Z(\mu,s)\end{array}\!\right)\, ds
\!&=&\! \frac{\mu}{|\alpha|}\ln(1/|\mu|)
\left\{\!
\left(\begin{array}{c}   
y_1(-\pi/2)+ y_1(\pi/2)\\
z_1(-\pi/2)+ z_1(\pi/2)
\end{array}\right)+o(1)\!\right\} \ .
\end{eqnarray} 
On the interval $[\pi/2,3\pi/2]$, the corresponding formulae hold with 
$\alpha$ replaced by $-\alpha$ on the right hand sides.
\end{lemma}

The appearance of the absolute values of $\mu$ and $\alpha$
plays a crucial role in the analysis of the first variation
of ${\cal L}$ in Proposition~\ref{prop:L1}.
We also need the following higher order expansion:

\begin{lemma}\label{lem:2ndorder}
Assume that a vector-valued function $\Xvec(s)$
on $[-\pi/2,\pi/2]$ satisfies
\begin{equation}\label{eq:expans2}
     \Xvec(\mu,s) =
\left(\begin{array}{c} \alpha\cos{(s)}\\ 0\\ 0\end{array}\right) 
+ \mu \xvec_1(s) + \mu^2 \xvec_2(s) + o(\mu^2)\end{equation}
in $H^1$, with $\xvec_1(-\pi/2)= \xvec_1(\pi/2)= 0$.  Then,
in the notation of Eq.~(\ref{eq:components-2}),
\begin{eqnarray}\label{eq:2ndx} 
\int_{-\pi/2}^{\pi/2}\frac{X(\mu,s)}
{|\Xvec(\mu,s)|}\, ds
&=&    
\pi\,\mbox{\rm sign}(\alpha)
+ \frac{\mu^2}{|\alpha|} \bigl(x_2(-\pi/2)+ x_2(\pi/2)\bigr)\\
&&\hskip -3cm - \frac{\mu^2}{\alpha}\bigl(|\xvec_2(-\pi/2)|+ 
|\xvec_2(\pi/2)|\bigr)
-\mbox{\rm sign}(\alpha)\frac{\mu^2\,}{2\alpha^2}\int_{-\pi/2}^{\pi/2}
\frac{1}{\cos^2 (s)} \bigl(y_1^2(s)+z_1^2(s)\bigr)\,ds \
   +o(\mu^2)\,,\nonumber
\end{eqnarray}
and
\begin{equation}\label{eq:2ndyz} 
 \int_{-\pi/2}^{\pi/2}\frac{1}{|\Xvec(\mu,s)|}\left(\begin{array}{l}
Y(\mu,s)\\Z(\mu,s)\end{array}\right)\,ds
= \frac{\mu}{|\alpha|}
\int_{-\pi/2}^{\pi/2}\frac{1}{|\cos (s)|} \left(\begin{array}{l}
y_1(s)\\z_1(s)\end{array}\right)\,ds+o(\mu)\,.
\end{equation}
On the interval $[\pi/2,3\pi/2]$, the corresponding formulae hold with 
$\alpha$ replaced by $-\alpha$ on the right hand sides.
\end{lemma}

\noindent {\bf Remark}: Since $\Xvec(s)$ is an $H^1$-function
with $\xvec_1(-\pi/2)=\xvec(\pi/2)=0$,
the integrals in Eq.~(\ref{eq:2ndx}) and Eq.~(\ref{eq:2ndyz})
are finite by Lemma~\ref{lem:lowerbound}.

\bigskip\noindent 
The proofs rely on the well-known fact that $H^1$-functions 
on the circle are bounded and
H\"older continuous with exponent $1/2$. We will need
the slightly stronger estimate
\begin{equation} \label{eq:Holder}
|x(s)-x(t)|
\leq \int_t^s |\frac{dx}{dt}(s')|\,ds'
\leq \left(\int_t^s ds'\right)^{1/2}\left(\int_t^s
     \left(\frac{dx}{dt}\right)^2(s')ds'\right)^{1/2}
= |s-t|^{1/2}o(1)\ .
\end{equation}
Since $ F(t)\equiv \int_0^t (\frac{dx}{dt})^2(s')\, ds'$
is uniformly continuous in $t$,
the $o(1)$ estimate holds uniformly in
$s$ and $t$.

\bigskip\noindent {\sc Proof of Lemma \ref{lem:1storder}.} \ 
Let $\Xvec(\mu,s)$ be of the form given in Eq.~(\ref{eq:xvecform1}),
and use the notation in Eq.~(\ref{eq:components-1})
for the component functions.
By the scaling invariance of the integrand, we
may replace $\alpha$ with 1 and $\mu$ with $\mu/\alpha$ without changing the 
values of the integrals. We also
assume that $\mu>0$, replacing $\mu$ with $-\mu$ and $\xvec_1$ 
with $-\xvec_1$ if necessary. 

Let us consider the resulting integral in the half-interval
$[0,\pi/2]$, beginning with a neighborhood of $\pi/2$
where the denominators are small.  For $s\in [\pi/2-\mu/\delta(\mu),
\pi/2]$ and with $\delta =\delta(\mu)= o(1)$ to be further specified
below, we see with the Taylor expansion of the cosine and the H\"older
continuity of the $H^1$-function $\xvec_1$ that
\begin{eqnarray*}
\Xvec(\mu,s) &=& \left(\begin{array}{c} \pi/2-s + O(s-\pi/2)^3\\ 0\\0
\end{array}\right) + \mu \bigl(\xvec_1(\pi/2) + o(s-\pi/2)^{1/2}\bigr)
+ o(\mu) \\ &=& \left(\begin{array}{c} \pi/2-s +\mu x_1(\pi/2)\\
\mu y_1(\pi/2)\\\mu z_1(\pi/2)
\end{array}\right) + O(\mu^3\delta^{-3})+ o(\mu^{3/2}\delta^{-1/2})+o(\mu)
\\
&=:& \vvec(\pi/2-s) + \left\{
O(\mu^3\delta^{-3})+o(\mu^{3/2}\delta^{-1/2})+o(\mu)\right\}\ .
\end{eqnarray*}
In the second step, we have used that $|s-\pi/2| \le \mu
/\delta$.   We may neglect contributions to
the integrals over the set
$$
\Delta = \Delta(\mu) :=
\bigl\{ s\in [0,\pi/2]: |\pi/2-s+\mu x_1(\pi/2)|\le \mu \delta(\mu)\bigr\}\ ,
$$
because the integrands are bounded, and the
measure of $\Delta$ is $o(\mu)$. On the complement of
$\Delta$ we use the inequality that for any pair of 
vectors $\vvec,\wvec$ with $|\vvec| \ge 2 |\wvec|>0$, 
\begin{equation} \label{eq:auxil}
\left\vert \frac{\vvec+\wvec}{|\vvec+\wvec|}-\frac{\vvec}{|\vvec|}\right|
\le 4\frac{|\wvec|}{|\vvec|}.
\end{equation}
We apply this to $\vvec(\pi/2-s)$ and $\wvec(\mu,s) =
O(\mu^3\delta^{-3})+o(\mu^{3/2}\delta^{-1/2})+ o(\mu)$ outside
of $\Delta$ with $\delta= \delta(\mu)$ now chosen so that
$\|\wvec\|_{\infty}/(\mu\delta^{2})= o(1)$, which is the case if
$\delta(\mu)$ exceeds $\mu^{1/5}$
and $o(\mu)/\mu\delta^2= o(1)$ where the $o(\mu)$-term refers to 
that in the expansion in Eq.~(\ref{eq:xvecform1}) and
$\delta(\mu)$ itself is still $o(1)$.  We obtain
\begin{eqnarray}
\lefteqn{ \int_{\pi/2-\mu/\delta(\mu)}^{\pi/2}
\frac{\Xvec(\mu,s)}{|\Xvec(\mu,s)|} - 
          \left(\begin{array}{c} 1 \\ 0 \\ 0 \end{array}\right)
\, ds }\nonumber\\ 
&=& \int_{[\pi/2-\mu/\delta(\mu),\pi/2]\setminus\Delta}
\left\{ \frac{\vvec(\pi/2-s)}{|\vvec(\pi/2-s)|}-
\left(\begin{array}{c} 1 \\ 0\\0 \end{array}\right) +\delta(\mu) o(1)\, 
    \right\}\,ds +o(\mu)  \nonumber \\
&=& \int_{0}^{\mu/\delta(\mu)}
\left\{ \frac{\vvec(s)}{|\vvec(s)|}
-\left(\begin{array}{c} 1 \\ 0 \\0\end{array}\right) \, \right\}ds+ o(\mu)\ .
\label{eq:x1y1}
\end{eqnarray}
The $x$-component of the integral in the last line of Eq.~(\ref{eq:x1y1})
is elementary
and equals
\begin{eqnarray}
&&\hskip -4cm
\int_{0}^{\mu/\delta(\mu)}
    \left\{\frac{s+ \mu x_1(\pi/2)}
    {\bigl(s+\mu x_1(\pi/2))^2+\mu^2y_1(\pi/2)^2+\mu^2z_1(\pi/2)^2\bigr)^{1/2}} -1  \right\}\, ds \nonumber\\
 &&= \sqrt{(s+\mu x_1(\pi/2))^2+\mu^2 y_1(\pi/2)^2+\mu^2 z_1(\pi/2)^2}-s
\big\vert_{s=0}^{\mu/\delta(\mu)}\nonumber \\
&&=  \mu \bigl(x_1(\pi/2)-|\xvec_1(\pi/2)|\bigr)+ o(\mu)\ .
\label{eq:x1-a}
\end{eqnarray}
The $y$-and $z$-components of the integral in Eq.~(\ref{eq:x1y1}) 
are computed similarly, e.g.,
\begin{eqnarray}
\label{eq:y1-a}
\lefteqn{\int_{0}^{\mu/\delta(\mu)}
\frac{ \mu y_1(\pi/2) }{((s+ \mu x_1(\pi/2))^2+\mu^2 y_1(\pi/2)^2+\mu^2 z_1(\pi/2)^2 )^{1/2}}
\, ds}  \nonumber \\
&=& \mu y_1(\pi/2)\ln{\left(s+ \mu x_1(\pi/2)+\left((s+ \mu x_1(\pi/2))^2+\mu^2 y_1(\pi/2)^2+\mu^2 z_1(\pi/2)^2\right)^{1/2} \right)\big|_0^{\mu/\delta(\mu)}}\nonumber\\
&=& \mu \ln{(\frac{1}{\delta(\mu)})}y_1(\pi/2)+O(\mu).
\end{eqnarray}
The error of order $O(\mu)$ reflects the shift 
of the zero in the denominator by $\mu x_1(\pi/2)$.
For the remaining part of the interval,
the cosine dominates the denominator, and one finds 
for the $x$-component that
 \begin{eqnarray}
&&\hskip -4cm \int_{0}^{\pi/2-\mu/\delta(/\mu)}
\Bigl\{ \frac{\cos(s)+ \mu x_1(s) + o(\mu)}
    {|\Xvec(s)|}-1\Bigr\} \, ds \nonumber\\
   &=& \int_0^{\pi/2-\mu/\delta(\mu)}
   O\left(\frac{\mu y_1(s)}{\cos(s)+\mu x_1(s)+o(\mu)}\right)^2\, ds
         \nonumber\\
     &=& O(\mu\delta(\mu))= o(\mu)\ .
\label{eq:x1-b}
   \end{eqnarray}
We have used that $\xvec_1(s)$ is uniformly bounded. 
For the $y$-component, we have
\begin{eqnarray}
&& \hskip -2cm \int_0^{\pi/2-\mu/\delta(\mu)}\frac{\mu y_1(s) + o(\mu)}
{|\Xvec(s)|}\, ds\nonumber\\
 &=&  \int_0^{\pi/2-\mu/\delta(\mu)} 
\frac{\mu y_1(\pi/2) + \mu o((\pi/2-s)^{1/2}) + o(\mu)}
{\cos(s)} \bigl(1+  o(1)\bigr)\, ds \nonumber\\
&=& \mu y_1(\pi/2)\ln\left(\sec(s)+\tan(s)\right)\big|_0^{\pi/2-\mu/\delta(\mu)}
  + o(\mu\ln(1/\mu))\nonumber\\
&=& -\mu \ln\bigl(\mu/\delta(\mu)\bigr)y_1(\pi/2) + o(\mu\ln(1/\mu)\ ,
\label{eq:y1-b}
\end{eqnarray}
where we have again exactly evaluated the integral
and expanded the result. The $z$-component is analyzed in the same way..

Adding Eqs.~(\ref{eq:x1-a}) and~(\ref{eq:y1-a}) to
Eqs.~(\ref{eq:x1-b}) and~(\ref{eq:y1-b}) respectively, we get that
$$
\int_0^{\pi/2} \left\{\frac{\Xvec(\mu,s)}{|\Xvec(\mu,s)|} - 
     \left(\begin{array}{c} 1 \\ 0 \\ 0 \end{array}\right) \right\}\,ds 
=  \left(\begin{array}{l}
  \mu\bigl\{x_1(\pi/2) -|\xvec_1(\pi/2)|\bigr\} + o(\mu)\\
 \mu \ln{(1/\mu)}y_1(\pi/2)  + o(\mu\ln(1/\mu))\\
 \mu \ln{(1/\mu)}z_1(\pi/2)  + o(\mu\ln(1/\mu))
\end{array}\right).
$$
To obtain Eq.~(\ref{eq:1storderx}), we repeat the computation
for the interval $[-\pi/2,0]$ and add the results. 
The claim for the interval $[\pi/2,3\pi/2]$ follows by replacing 
$\Xvec(s)$ with $-\Xvec(s-\pi)$.
\hfill\qed

\bigskip 
\noindent {\sc Proof of Lemma \ref{lem:2ndorder}.}\ 
Here, we assume that $\Xvec(\mu,s)$ has the expansion 
in Eq.~(\ref{eq:expans2}) and $\xvec_1(-\pi/2)= \xvec(\pi/2)=0$. 
We may assume by scaling that $\alpha=1$ and $\mu>0$.
Let us use again the notation in
Eq.~(\ref{eq:components-2}) to denote
the components of the various vector-valued functions.

We will expand the integrand and partition the
interval of integration as in the proof of Lemma~\ref{lem:1storder}.
By Eq.~(\ref{eq:Holder}), there is a function $m(s)= o(s^{1/2})$
such that $|\xvec_i(s)-\xvec_i(t)| \leq m(|s-t|)$, for $i=1,2$.  Let
$\delta= \delta(\mu)= o(1)$ to be further specified below.  On
$[\pi/2-\mu^2/\delta,\pi/2]$ we expand
\begin{eqnarray}\label{eqmess}
\Xvec(\mu,s) &=& \left(\begin{array}{c}\pi/2-s +
O(s-\pi/2)^3\\0\\0\end{array}\right) + \mu^2\xvec_2(\pi/2) + \mu
O(m(s-\pi/2)) + o(\mu^2)\nonumber \\ 
&=& \left(\begin{array}{c}
\pi/2-s+ \mu^2 x_2(\pi/2) \\ \mu^2y_2(\pi/2)\\\mu^2z_2(\pi/2)
\end{array}\right)
+ O(\mu^6/\delta^3) + \mu m(\mu^2/\delta) + o(\mu^2)\nonumber
\\
 &=:&
\vvec(\pi/2-s) + \{O(\mu^6/\delta^3) + \mu m(\mu^2/\delta) +
o(\mu^2)\}.
\end{eqnarray}
At this point we choose $\delta = \delta(\mu)$ so that
$\mu^4\delta^{-5}= o(1)$,\, $m(\mu^2/\delta)\mu^{-1}\delta^{-2}=o(1)$
and that $\mu^{-2}\delta^{-2}o(\mu^2)= o(1)$, still
keeping $\delta(\mu)= o(1)$.  This will ensure that the sum 
of the last three terms of Eq.~(\ref{eqmess}) divided by 
$|\vvec(\pi-s)|$, is no bigger than $\delta(\mu)\times o(1)$ outside of
 $\Delta$ defined by
$$
\Delta=\Delta(\mu) = \left\{s\in [0,\pi/2]: |s-\pi/2 + \mu^2 x_2(\pi/2)|
\le \delta\mu^2\right\}\ .
$$
We  again neglect the integral over $\Delta$, since
$$
\left| \int_\Delta 
\left\{ \frac{\Xvec(\mu,s)}{|\Xvec(\mu,s)|} 
            - \left(\begin{array}{c} 1 \\ 0 \\ 0 \end{array}\right)
\right\}
\, ds\right|
    \le 8 \mu^2\delta\ .
$$ 
We also apply the vector inequality Eq.~(\ref{eq:auxil}) again; we obtain
\begin{eqnarray*}
&& \hskip -2cm \int_{\pi/2-\mu^2/\delta}^{\pi/2}
\left\{
\frac{\Xvec(\mu,s)}{|\Xvec(\mu,s)|}- 
\left(\begin{array}{c}1 \\ 0 \\ 0 \end{array}\right) \right\}\, ds\\
&=& \int_{[\pi/2-\mu^2/\delta,\pi/2]\setminus\Delta}
\left\{ \frac{\vvec(\pi/2-s)}{|\vvec(\pi/2-s)|} -
\left(\begin{array}{c}1\\0\\0\end{array}\right) \right\}
\, ds + o(\mu^2).
\end{eqnarray*}
where the last $o(\mu^2)$-term is simply $\delta(\mu)\times o(1)\times
\mu^2/\delta(\mu)$ coming from the integral of the vector inequality,
and from neglecting the integral over $\Delta$. The integral on the
right side of this last expression is done explicitly and then
estimated as in the proof of the previous lemma, giving
\begin{equation}\label{eq:mint}
\int_{\pi/2-\mu^2/\delta}^{\pi/2}
\frac{\Xvec(\mu,s)}{|\Xvec(\mu,s)|}- 
\left(\begin{array}{c} 1 \\ 0 \\ 0 \end{array}\right) \, ds
= \left(\begin{array}{l}
   \mu^2\bigl(x_2(\pi/2)-|\xvec_2(\pi/2)|\bigr)+o(\mu^2)\\
    \mu^2\ln(1/\delta) y_2(\pi/2)+O(\mu^2)\\
    \mu^2\ln(1/\delta) z_2(\pi/2)+O(\mu^2)
\end{array}\right).
\end{equation}

When $s\in [0,\pi/2-\mu^2/\delta]$,
the cosine dominates both the numerator and denominator,
\begin{eqnarray*}
\frac{\mu\xvec_1(s)+\mu^2\xvec_2(s)+o(\mu^2)}{\cos(s)}
&=& \frac{\mu(\xvec_1(s)-\xvec_1(\pi/2))+\mu^2\xvec_2(s)+o(\mu^2)}
{\cos(s)}\nonumber\\
&=& O(\delta^{1/2}) = o(1).
\end{eqnarray*}
For the $x$-component of the integral, we have
\begin{eqnarray}\label{eqnlaintx}
\lefteqn{\int_0^{\pi/2-\mu^2/\delta}
\left\{\frac{X(\mu,s)}{|\Xvec(\mu,s)|}-1\right\}\,ds}\nonumber\\
&=& -\frac{1}{2}\int_0^{\pi/2-\mu^2/\delta} 
\frac{(\mu y_1(s)+O(\mu^2))^2+(\mu z_1(s)+O(\mu^2))^2}
   {\cos^2(s)}(1+o(1))\,ds\nonumber\\
&=& -\frac{\mu^2}{2}\int_0^{\pi/2}\frac{y_1^2(s)+z_1^2(s)}{\cos^2(s)}\, ds +o(\mu^2). 
\end{eqnarray}
In the last line we used Lemma~\ref{lem:lowerbound}
to see that
$y_1(s)/\cos(s)$ and $z_1(s)/\cos(s)$ are square integrable over the 
{\em entire interval} $[0,\pi/2]$,
so that extending the interval 
of integration introduces only an additional 
$\mu^2\times o(1)= o(\mu^2)$ error.
For the $y$-component of the integral, we get that
\begin{eqnarray}\label{eqnlainty}
 \int_0^{\pi/2-\mu^2/\delta}\frac{Y(\mu,s)}{|\Xvec(\mu,s)|}\,ds
&=& \mu\int_0^{\pi/2-\mu^2/\delta}\frac{y_1(s)+O(\mu)}{\cos(s)}(1+o(1))\,ds\nonumber\\
&=& \mu\int_0^{\pi/2}\frac{y_1(s)}{\cos(s)}\,ds +o(\mu)
\end{eqnarray}
and a similar expression for the $z$-component, where again extension
of the interval of integration introduces only an $o(\mu)$ error.
Collecting the results of
Eqs.~(\ref{eq:mint})-(\ref{eqnlainty}), we obtain
$$
   \int_{0}^{\pi/2}
\left(\frac{X(\mu,s)}{|\Xvec(\mu,s)|}-1\right)\,ds
 = \mu^2\left\{x_2(\pi/2)-|\xvec_2(\pi/2)|
-\frac{1}{2}\int_{0}^{\pi/2}\frac{y_1^2(s)}{\cos^2(s)}\,ds
\right\} +o(\mu^2)
$$
and
\begin{eqnarray*}
   \int_0^{\pi/2}\frac{1}{|\Xvec(\mu,s)|}\left(\begin{array}{l}
Y(\mu,s)\\Z(\mu,s)\end{array}\right)\,ds &=&  \mu \int_{0}^{\pi/2} 
         \frac{1}{\cos(s)}\left(\begin{array}{l}
y_1(s)\\z_1(s)\end{array}\right)\,ds +
o(\mu).
\end{eqnarray*}
To arrive at Eqs.~(\ref{eq:2ndx}) and (\ref{eq:2ndyz}), we repeat
the computations on the interval $[-\pi/2,0]$ and add the results.
The claim for the interval $[\pi/2,3\pi/2]$ follows by replacing
$\Xvec(s)$ with $-\Xvec(s-\pi)$.
\hfill\qed       

\newpage

%***********************************
\appendix

\section*{Appendix}

\section{Eigenvalues of a Sturm-Liouville operator}
\setcounter{equation}{0}
\setcounter{theorem}{0}

We provide an overview of the spectral theory for the operator 
$$
K_g = -\frac{d^2}{ds^2}+g\sec^2(s)$$
on $[-\pi/2,\pi/2]$, with Dirichlet 
boundary conditions at the endpoints
(cf. Methods of Theoretical Physics~\cite{MF}, P.M. Morse and 
H. Feshbach, Part I, p.388 and the discussion there of 
hypergeometric functions.)
Here, $g$ is a constant.  We first show that $K_g$ is 
bounded below for $g \ge -\frac{1}{4}$.

\begin{lemma}
\label{lem:lowerbound}
 Suppose that $w(t)$ is an $H^1$ function on
 $[0, a]$, vanishing at $t=0$ and $t=a>0$.  Then
$$
        \frac{1}{4}\int_0^a \frac{w(s)^2}{s^2}ds\leq \int_0^a
        \left(w'(s)\right)^2\,ds.
$$
\end{lemma}

\noindent {\sc Proof.}\  By scale invariance, it suffices
to consider the case $a=1$. We have that
\begin{eqnarray}
  0&\leq & \int_t^1
  \left(w'(s)-\frac{w(s)}{2s}\right)^2\,ds\nonumber\\
\nonumber   &=& \int_t^1 \left(\frac{dw}{ds}(s)\right)^2\,ds
    - \frac{1}{2}\int_t^1 \frac{d/ds\, w^2(s)}{s} \, ds 
  +\frac{1}{4}\int_t^1\left(\frac{w(s)}{s}\right)^2\,ds\ .
\end{eqnarray}
Integrating by parts in the second integral
and collecting terms, we  get
$$
\frac{1}{4}\int_t^1\left(\frac{w(s)}{s}\right)^2\,ds
\leq  \int_t^1 \left(w'(s)\right)^2\,ds 
-\frac{w^2(s)}{s}\Big\vert^{s=1}_{s=t}\ .
$$
By assumption, $w(1)=0$, and by Eq.~(\ref{eq:Holder}),
$w(t)=o(t^{1/2})$.
The desired conclusion follows by taking $t\to 0$.
\hfill \qed

The lemma implies that $K_g$ is bounded below for $g\ge -1/4$, because
$$ 
\inf_{s\in [0,2\pi]}
\left\{
\frac{1}{(\pi/2-s)^2} + \frac{1}{(3\pi/2-s)^2} -\sec^2(s) \right\}
>-\infty\ .  $$ 
Furthermore, $K_g$ has compact resolvent when $g>-1/4$,
since $K_g\ge -c_1(g)  d^2/ds^2 - c_2(g)I$ for some constants 
$c_1(g), c_2(g)>0$, and the positive operator $-d^2/ds^2$ 
has compact resolvent. Consequently,
the spectrum of $K_g$ consists 
of a nondecreasing sequence of eigenvalues $\lambda_0<\lambda_1\le \dots$
with $\lambda_n\to\infty$.  The ground state $\lambda_0$ is simple
by a Perron-Frobenius argument.

To solve the eigenvalue-eigenvector equation 
$$ K_gw(s)= \lambda w(s)\ ,$$
one can write $w= \cos^a(s)\phi(s)$ with 
\begin{equation} \label{eq:a-g}
a= \frac{1}{2}\left(1+\sqrt{1+4g}\right)
\end{equation}
and obtain a second order
differential equation for $\phi$. A 
substitution $ \xi = (1+\sin(s))/2$
results in the {\em hypergeometric} equation for $\phi$ regarded now
with a slight abuse of notation as a function of $\xi$ 
$$
 -\xi(1-\xi)\frac{d^2\phi(\xi)}{d\xi^2}
  +2(a+\frac{1}{2})(\xi-\frac{1}{2})\frac{d\phi(\xi)}{d\xi}
  +(a^2-\lambda)\phi(\xi)= 0.\\
$$
Expanding $\phi$ in a power series about $\xi= 0$, one obtains a
      hypergeometric series, 
$$ \phi(\xi)= \sum_{n=0}^{\infty}
      b_n \xi^n
$$
with the coefficients $b_n$ satisfying a two-term recursion relation,
$$
           b_{n+1}= \frac{(n+a)^2-\lambda}{(n+a+\frac{1}{2})(n+1)}b_n;
$$
 (The indicial equation gives that the series indeed should begin with
the $n=0$ term. The other solution leads to a function $w$ which is
not locally $H^1$ at $-\pi/2$, i.e., $\frac{dw}{ds}$ is not locally
square-integrable there). One finds that 
$$
      b_n= \frac{\Gamma(a+\frac{1}{2})}     {\Gamma(r_1)\Gamma(r_2)}\times 
           \frac{\Gamma(r_1+n)\Gamma(r_2+n)}{\Gamma(a+n+\frac{1}{2})\,n!},
$$
where$-r_1$ and $-r_2$ are the roots of the equation
$n^2+2an +a^2-\lambda= 0$.  Via Stirling's approximation, one can
infer from the expression for the $b_n$'s that $b_n \sim
n^{a-3/2}(1+{\cal O}(1/n))$ for $n$ large further implying that
$\phi(\xi)\sim (1-\xi)^{1/2-a}$ or that $w(s)$ would not be locally
square integrable in a neighborhood of $s=\pi/2$. (Alternatively 
this conclusion can be arrived at through well-known integral 
representations for hypergeometric functions.)  Thus $b_n$ must be 
eventually zero.  It follows from the recursion relation
that the eigenvalues $\lambda_n$ satisfy the quantization condition 
\begin{equation}
\label{eq:lambda-gegenbauer}
\lambda_n= (n+a)^2\ , \quad n= 0,1,.. 
\end{equation}
In particular, the ground state satisfies $\lambda_0=a^2 \ge 1/4$
for all $g>-1/4$.

The function $\phi_n(\xi)$ corresponding to $\lambda_n$ is a 
polynomial of degree $n$. In fact,
with the further transformation $z= 2\xi-1$, the equation for $\phi_n$
as a function of $z$ is that of a Gegenbauer polynomial,
$$
   (z^2-1)\frac{d^2}{dz^2}\phi_n(z)+(2a+1)z\frac{d}{dz}\phi_n(z)-
(2an+n^2)\phi_n(z)= 0
$$
with solution $\phi_n(z)= T_n^{a-\frac{1}{2}}(z)$, with well-known
orthogonality and normalization properties.  The resulting functions
$\{w_n(s)= \cos^a(s)T_n^{a-\frac{1}{2}}(\sin(s))\}$ are complete.
\hfill \qed

\bigskip \noindent {\bf Remark}: (1) 
Recalling the relationship between
the parameters $a$ and $g$ from Eq.~(\ref{eq:a-g}), we see that
Eq.~(\ref{eq:lambda-gegenbauer}) implies the lower bounds
\begin{equation} \label{eq:lambda01}
\left\{
\begin{array}{ll}
\lambda_0>\frac{1}{4}, \lambda_1 > 1\ ,
\quad & g>-\frac{1}{4}\\
\lambda_0 > 1\ , & g> 0\ .
\end{array}\right.
\end{equation}

(2) When $g\leq -1/4$,
the function $\cos^a(s)\phi(\frac{1+\sin(s)}{2})$ appearing in the
change of variables is no longer locally in $H^1$ and the above construction 
of the eigenfunctions and eigenvalues does not
apply.  For $g=-1/4$ we have
the sharp inequality
$$
 \frac{1}{4}\int_{-\pi/2}^{\pi/2}\sec^2(s)w^2(s)\,ds
\leq \int_{-\pi/2}^{\pi/2}\left(\frac{dw(s)}{ds}\right)^2\,ds
- \frac{1}{4} \int_{-\pi/2}^{\pi/2} w^2(s)\, ds
$$
for functions $w$ satisfying Dirichlet conditions at $\pm \pi/2$: Our
above analysis gives this result with the $1/4$ on the left side
replaced by $-g<1/4$, and taking $g\downarrow -1/4$ completes 
the argument.

%************************************* 

\end{document}